\documentclass[11pt]{article}

\usepackage[english]{babel}
\usepackage{graphicx}
\usepackage{subcaption}
\usepackage{amsmath,amssymb,amsthm}
\usepackage{siunitx}
\usepackage[margin=2.5cm]{geometry}
\usepackage{xcolor}

\newcommand{\R}{\mathbb{R}}

\newcommand{\C}{\mathbb{C}}
 \newcommand{\ds}{\displaystyle}

\usepackage[
  colorlinks = true,
  linkcolor  = blue,
  urlcolor   = blue,
  citecolor  = blue
]{hyperref}

\title{The Unified Transform for Burgers' Equation: Application to Unsaturated Flow in a Finite Interval}

\author{
Konstantinos Kalimeris$^{1}$, 
Leonidas Mindrinos$^{2}$ and
Athanasios Paraskevopoulos$^{1,2}$%
\thanks{Corresponding author: \texttt{at.paraskevopoulos@aua.gr}}
}

\date{
$^{1}$ Mathematics Research Center, Academy of Athens, Greece\\
$^{2}$ Department of Natural Resources Development and Agricultural Engineering,\\
Agricultural University of Athens, Greece
}

\begin{document}

\maketitle

\begin{abstract}
In this paper, we focus on one-dimensional vertical infiltration, assuming constant diffusivity and a quadratic relationship between hydraulic conductivity and water content. Under these assumptions, Richards’ equation reduces to Burgers’ equation, which we then linearize via the Hopf–Cole transformation. This turns the initial boundary value problem into a diffusion equation on a finite interval with mixed boundary conditions. To solve it, we use the Unified Transform Method (also known as the Fokas method). This approach gives an explicit integral representation of the solution, and when evaluated numerically, the results match classical Fourier series solutions exactly, but with better convergence and stability. Two examples from hydrological applications are examined.
\end{abstract}

\section{Introduction}

The problem of water flow in unsaturated porous media has attracted much attention since the early years of soil physics research, as it presents significant theoretical and experimental challenges. From a mathematical perspective, great effort has been devoted to the derivation of both analytical and numerical solutions that can be compared with experimental observations to enhance understanding of the underlying physical processes.

In this work, we focus on the vertical infiltration problem (one-dimensional) in an unsaturated porous medium. This process is described by Richards' equation \cite{richards}, which is a non-linear partial differential equation (PDE) that combines diffusion and transport phenomena. Since its introduction in the early 20th century, Richards' equation has been widely used to model various physical conditions. Despite its strong physical foundation, the nonlinear nature of the equation significantly complicates the derivation of analytical solutions.

A common approach to simplify this model is to impose specific functional forms for the soil hydraulic parameters. Typically, the soil water diffusivity ($D$) is assumed constant and the hydraulic conductivity $(K)$ is a function of the water content $(\theta)$. Variations also arise depending on the imposed boundary conditions, which may represent flooding or rainfall scenarios and are commonly expressed as Dirichlet or Robin conditions, respectively. Additional restrictions are made if the soil domain is treated as bounded or semi-infinite. The following cases admit approximate analytical solutions for both flooding and rainfall conditions: 
\begin{description}
    \item[Case 1] $K(\theta) =$ constant \cite{Bra73, CarJae59}.
     \item[Case 2] $K(\theta) =$ linear function of $\theta$ \cite{Bra73, philip66}.
      \item[Case 3] $K(\theta) =$ quadratic function of $\theta$ \cite{clothier, philip1974}.
\end{description}

The first case corresponds to neglecting gravitational effects and reduces Richards’ equation to the linear diffusion equation. Case 2 simplifies the governing equation to the linear advection–dispersion equation, whereas in the third case, the PDE reduces to a form of Burgers’ equation. The first two cases have been widely examined in the literature; here, we focus on the third case. Although the resulting equation remains nonlinear, it can be transformed into a linear diffusion equation through the Hopf–Cole transformation \cite{hopf50}. Consequently, the original problem can be reformulated as an initial–boundary value problem (IBVP) for the heat equation, similarly to Case 1 \cite{hills, philip1974}.

A variety of analytical methods for the heat equation have been developed, among which separation of variables and Fourier series expansions are the most widely used. These classical techniques have been successfully applied to both infiltration and drainage problems in finite or semi-infinite domains under various boundary conditions \cite{hills, warrick1995analytical, basha2002burgers}. Indeed, all previously reported analytical studies rely on such eigenfunction-based representations.

Despite their effectiveness, these approaches impose restrictive assumptions on problem formulation. In particular, they are typically limited to finite spatial intervals and require homogeneous or simple local boundary conditions. As a result, the approximate analytical solution often leads to lengthy derivations and series solutions that may converge slowly, especially for short times or near boundaries where compatibility conditions might not apply. These limitations, reflecting the lack of uniform convergence of the solution to their traces on the boundary of the IBVP, motivate the development of alternative analytical frameworks that can address such problems in a more direct and unified manner.

In this study, we apply the Unified Transform Method which was introduced in the seminal work of Fokas in 1997 \cite{fokas97}; hence, the method is now known as the Fokas method. This approach yields an explicit integral representation of the solution relying on the naturally arosen spectral transforms of the initial and bundary conditions. We refer to the book \cite{FK22} and the recent work \cite{chatziafratis2025fokas} where the Fokas method is applied to various evolution PDEs and the insight of the appearance of continuous spectrum is developed for finite domains, in contrast to the classically derived discrete Fourier spectrum. As a result, the analytical formulation is significantly simplified, and the resulting integral representation is often better suited for numerical evaluation, exhibiting improved stability and convergence properties when compared with truncated Fourier series solutions.

Case~1 on the half-line was first studied in the original work \cite{fokas97}, followed by its extension to finite intervals \cite{FokPel05}. Further applications have been reported in later studies \cite{cha22, KO20, KOD24, man13, KalMin25b}. The linear evolution equation (see Case 2) has also been extensively examined on the half-line \cite{barros19, fokas02a, fokas02b} and on the finite interval \cite{arg24, hwang24, KalMin25} under various boundary conditions. To our knowledge, the application of the Fokas method to Burgers' equation (Case 3) has so far been considered only on the half-line \cite{foklil}. The present work addresses this gap by extending the method to a finite interval with mixed boundary conditions. The possibility of investigating inverse problems in the vein of the above-mentioned literature, namely the determination of boundary controllers associated with physically relevant problems, will be considered in future work.

The structure of the paper is as follows. In \autoref{sec_formulation}, we derive Burgers’ equation for the specific physical problem and reduce it to a diffusion equation on a finite interval under mixed boundary conditions. The resulting IBVP is then solved analytically using the Unified Transform Method. The numerical implementation of the integral representation is presented in \autoref{sec_implementation}. In the final section, the analytical solution is compared with existing Fourier series solutions under physically relevant conditions, demonstrating agreement while highlighting the advantages of the proposed approach.

\section{Problem formulation}\label{sec_formulation}

In soil physics, the vertical water flow in an unsaturated porous medium is modelled through the one-dimensional non-linear Richards' equation \cite{richards}
\begin{equation}\label{richards}
\frac{\partial \theta}{\partial t} = \frac{\partial}{\partial x} \left(D(\theta) \frac{\partial \theta}{\partial x} \right) - \frac{\partial K (\theta)}{\partial x},
\end{equation}
where $\theta$ is the water content, $K$ is the hydraulic conductivity and $D$ describes water diffusivity.

There are several ways to describe how the material parameters $D$ and $K$ depend on $\theta.$ In this work, we focus on Case 3 of the Introduction, by assuming constant diffusivity ($D>0$) and that $K$ is a quadratic function of $\theta$ \cite{philip1974}:
\begin{equation}\label{condu}
K (\theta) = a (\theta + b)^2.
\end{equation} 
Then, equation \eqref{richards} takes the form
\begin{equation}\label{richards2}
\frac{\partial \theta}{\partial t} = D \frac{\partial^2 \theta}{\partial x^2}  - 2a( \theta +b ) \frac{\partial \theta}{\partial x}.
\end{equation}
 This is the well-known Burgers' equation \cite{burgers}, which is still non-linear but can be linearized using the Hopf-Cole transformation. 

Following \cite{hills} one can show that \eqref{richards2} reduces to the linear diffusion equation
\begin{equation}\label{heat}
\frac{\partial w}{\partial t} = D \frac{\partial^2 w}{\partial x^2},
\end{equation}
for the new variable
\begin{equation}\label{new_var}
w = e^{-\ds\frac{a}{D}u}, \quad \mbox{where} \quad \frac{\partial u}{\partial x} = \theta + b. 
\end{equation}

Thus, once \eqref{heat} is solved for $w$ (subject to appropriate initial and boundary conditions), the water content is given by
\begin{equation}\label{solution}
\theta  = -\frac{D}{a} \frac{\partial}{\partial x} (\ln w) - b =   -\frac{D}{a} \frac{w_x}{w}  - b.
\end{equation}

We consider \eqref{richards2} in a finite interval $x\in[0,L]$ together with  a constant initial condition of the form
\begin{equation}\label{initial}
\theta (x,0) = \theta_0 , \quad 0<x<L.
\end{equation}

At the surface $x=0,$ we assume given time-dependent flux $q$, modelled through Darcy's law
\begin{equation}\label{darcy}
K(\theta) - D (\theta) \frac{\partial \theta}{\partial x} (0,t) = q(t), \quad t>0. 
\end{equation}
and at the bottom ($x=L$) we  impose a Dirichlet boundary condition (given moisture)
\begin{equation}\label{bc1}
\theta (L,t) = \theta_L , \quad t>0,
\end{equation}
with $\theta_0 \geq \theta_L >0.$

Equation \eqref{darcy} using \eqref{condu} takes the form of a non-linear Robin-type boundary condition
\begin{equation}\label{darcy2}
a (\theta(0,t) + b)^2 - D  \frac{\partial \theta}{\partial x} (0,t) = q(t), \quad t>0. 
\end{equation}

The initial \eqref{initial} and boundary conditions \eqref{bc1} and \eqref{darcy2}  have to be transformed using the new variables \eqref{new_var}. Thus, the final form of the IBVP is given by \cite{hills}
\begin{subequations}\label{ibvp}
\begin{alignat}{3}
\frac{\partial w}{\partial t}  &= D \frac{\partial^2 w}{\partial x^2},  \quad && 0<x <L, \, t>0,  \label{ibvp1}\\
w (x,0) &= w_0 (x), \quad && 0 <x <L, \label{ibvp2}\\ 
w (0,t) &= f(t), \quad &&t>0, \label{ibvp3}\\
\frac{\partial w}{\partial x} (L,t) + C\, w (L,t) &= 0, \quad &&t>0, \label{ibvp4}
\end{alignat}
\end{subequations}
where 
\begin{equation}\label{ib_functions}
w_0 (x) = e^{-\ds\frac{a}{D}(\theta_0 + b)x}, \qquad
f(t) = e^{\ds\frac{a}{D} \int_0^t q(s) ds},
\end{equation}
and $C = \ds\frac{a}{D} (\theta_L + b).$

\subsection{Analytical solution of the IBVP \eqref{ibvp}}

We are interested in finding the analytical solution of \eqref{ibvp} using the Fokas method. This is an IBVP for the heat equation in a finite interval with mixed boundary conditions.  In \cite{KalMin25}, two of the authors derived the solution of the vertical infiltration problem considering linear soil and various boundary conditions. The corresponding IBVP had the form
\begin{subequations}\label{bvp}
\begin{alignat}{3}
\frac{\partial \theta}{\partial t} + K_0 \frac{\partial \theta}{\partial x} &= D_0 \frac{\partial^2 \theta}{\partial x^2},  \quad && 0<x <L, \, t>0,  \label{bvp1}\\
\theta (x,0) &= \theta_0 (x), \quad && 0 <x <L, \label{bvp2}\\ 
\theta (0,t) - \alpha \frac{\partial \theta}{\partial x} (0,t) &= f(t),  \quad &&t>0, \label{bvp3}\\
\theta (L,t) -\beta \frac{\partial \theta}{\partial x} (L,t) &= g(t),  \quad &&t>0. \label{bvp4}
\end{alignat}
\end{subequations}

To obtain the analytical solution of \eqref{ibvp} we set $K_0 = 0, \, D_0 = D,$ in \eqref{bvp1}, $\alpha = 0$ in \eqref{bvp3} and $\beta = -\frac{1}{C},$ $g(t)=0$ in \eqref{bvp4}. Thus, we get (see \cite[Eq. (16)]{KalMin25})
\begin{equation}\label{solution_int}
\begin{aligned}
w (x, t)&=\frac{1}{2 \pi} \int_{\R} e^{i \lambda x-D \lambda^2 t} \hat{w}_0(\lambda) d \lambda \\
 & \phantom{=}-\frac{1}{2\pi} \int_{\partial D_{+}} \frac{e^{-D \lambda^2 t}}{ \Delta (\lambda, -L)}\left[ \sin (\lambda x) \left( i \lambda +C \right) e^{i  \lambda L} \hat{w}_0(\lambda)  \right. \\
&\phantom{=}\left. +\Delta (\lambda, x-L)  \hat{w}_0(-\lambda) + i 2 \lambda D \,  \Delta (\lambda, x-L) \tilde{f} (D\lambda^2, t) \right] d \lambda,
\end{aligned}
\end{equation}
where 
\begin{equation}\label{def:Delta}
\Delta (\lambda, y)=  \lambda   \cos \left( \lambda y \right) - C\sin \left( \lambda y \right)
\end{equation}
and the integral transforms appearing in \eqref{solution_int} are defined by
\begin{equation}\label{eq_transforms}
\hat{w}_0(\lambda) = \int_0^L e^{-i \lambda x} w_0 (x)dx, \qquad
\tilde{f} (k, t)  = \int_0^t e^{k s} f(s) ds, 
\end{equation}
namely the Fourier transform of the initial condition and the $t-$transform of the boundary conditions, respectively. 

The complex contour $\partial D_{+}$ is the positively oriented boundary of the domain $D_{+} = \big\{ \lambda \in \C: \operatorname{Re}(\lambda^2) < 0, \ \operatorname{Im}(\lambda) > 0 \big\}.$ There is freedom in choosing the boundary of $D_+$ to simplify the computation of the corresponding integral by avoiding contributions from any poles; see the relative discussion in Remark 2.1 in \cite{KalMin25}. Different cases will be examined in the next section. Alternatively, one can compute \eqref{solution_int} directly using numerical integration; however, this approach can be time-consuming and may require special treatment at early times or near the boundaries.

\section{Numerical implementation}\label{sec_implementation}

In this section, we deal with the numerical evaluation of the integral representation of the solution \eqref{solution_int}.

Without being restrictive for the realization of the numerical scheme, and considering a physically relevant scenario, we assume constant flux $q(t)=q$ and the formulas in \eqref{ib_functions} are simplified as
\[
w_0 (x) = e^{- A x}, \qquad
f(t) = e^{ B t},
\]
with 
\[
A = \frac{a}{D}(\theta_0 + b), \qquad B= \frac{a}{D}q.
\]

Then, from the definitions \eqref{eq_transforms} we obtain
\begin{equation}
\begin{aligned}
\hat{w}_0(\lambda)
&=\int_0^L e^{-A x}e^{-i\lambda x}\,dx
=\frac{1-e^{-(A +i\lambda)L}}{A +i\lambda}, \\
\tilde{f}(D\lambda^2 ,t)
&=\int_0^t e^{D\lambda^2 s}e^{Bs}\,ds
=\frac{e^{(D\lambda^2+B)t}-1}{D\lambda^2+B}.        
\end{aligned}
    \end{equation}
and substituting in \eqref{solution_int} we get the form
\begin{equation}\label{solution_integ}
\begin{aligned}
     w(x,t)&=\frac{1}{2\pi} \int_{-\infty}^{\infty} e^{i\lambda x- D\lambda^2 t} \frac{1-e^{-(A +i\lambda)L}}{A +i\lambda} \: d\lambda \\
    &\phantom{=}- \frac{1}{2\pi} \int_{\partial D^{+}} \frac{e^{-D\lambda^2 t}}{\Delta (\lambda, -L)}\left[ e^{i\lambda L} (i\lambda+C) \sin{(\lambda x)}\, \frac{1-e^{-(A +i\lambda)L}}{A+i\lambda} \right.\\
    &\phantom{=}\left. + \Delta (\lambda, x-L) \, \frac{1-e^{-(A -i\lambda)L}}{A -i\lambda}+ 2i \lambda D \, \Delta (\lambda, x-L) \frac{e^{(D\lambda^2+B)t}-1}{D\lambda^2+B}\right] \,d\lambda,
\end{aligned}
\end{equation}
where $\Delta(\lambda,y)$ is defined in \eqref{def:Delta}.

The analyticity and boundedness of the integrand 
of the first integral on the right-hand side of the above equation 
allows the deformation from the real line to the contour $\partial D^{+}$ in the following way
\begin{align*}
 \int_{-\infty}^{\infty} e^{i\lambda x- D\lambda^2 t} \frac{1-e^{-(A +i\lambda)L}}{A +i\lambda}  d\lambda &=  \int_{\partial D^+} e^{i\lambda x- D\lambda^2 t} \frac{1}{A +i\lambda}  d\lambda \\
&\phantom{=}- \int_{\partial D^-} e^{i\lambda x- D\lambda^2 t} \frac{e^{-(A + i \lambda)L}}{A +i\lambda}  d\lambda
\\    
&= \int_{\partial D^+} e^{- D\lambda^2 t} \Bigg[\frac{e^{i \lambda x}}{A +i\lambda} -   \frac{e^{-A L}e^{i\lambda (L-x)}}{A -i\lambda} \Bigg] \: d\lambda, 
\end{align*}
where inn the last step we used the change of variables $\lambda \rightarrow -\lambda$, to map $\partial D^-$ to $\partial D^+$. Combining all the above and by deforming $\partial D^+$ to $C^+,$ where $C^+$ is any simply connected contour with two asymptotes between the real axis and $\partial D^+$ (see the relative discussion in Remark 2.1 in \cite{KalMin25}),  we end up with the formula
\begin{equation}\label{eq_final}
    w(x,t) = w_1 (x,t) + w_2 (x,t),
\end{equation}
where
\begin{align*}
    w_1 (x,t) &=  \frac{1}{2 \pi} \int_{C^+ } \frac{e^{-D \lambda^2t}}{\Delta (\lambda, -L)}\left[\frac{e^{i \lambda x}(A-i \lambda)-e^{-A L}e^{i \lambda (L-x)}(A +i \lambda)}{A^2+\lambda^2}\Delta (\lambda, -L)\right.\\
    &\phantom{=}- e^{i \lambda L}(i \lambda+C) \sin{(\lambda x)}\frac{1-e^{- A L}e^{-i\lambda L}}{ A+i \lambda}\\
    &\phantom{=}\left.- \Delta(\lambda, x-L) \frac{1-e^{- A L}e^{i \lambda L}}{A- i \lambda}+2 i \lambda D \, \Delta(\lambda, x-L) \frac{1}{D\lambda^2+B}\right] d \lambda, \\
    w_2 (x,t) &= -\frac{1}{2\pi}\int_{C^{+}} \frac{1}{\Delta (\lambda, -L)} \left[ 2 i \lambda D \, \Delta(\lambda, x-L) \frac{e^{Bt}}{D\lambda^2+B}\right] d \lambda .
\end{align*}

After computing $w$ using \eqref{eq_final}, we substitute it into \eqref{solution} to obtain $\theta$. The spatial derivative can be formally evaluated, resulting in
\begin{align*}
   \frac{\partial w_1}{\partial x}  (x,t) &=  \frac{1}{2 \pi} \int_{C^+ } \frac{e^{-D \lambda^2t}}{\Delta (\lambda, -L)}\left[i \lambda \frac{ e^{i \lambda x}(A-i \lambda)+e^{-A L}e^{i \lambda (L-x)}(A +i \lambda)}{A^2+\lambda^2}\Delta (\lambda, -L)\right.\\
    &\phantom{=}- e^{i \lambda L} \lambda (i \lambda+C) \cos{(\lambda x)}\frac{1-e^{- A L}e^{-i\lambda L}}{ A+i \lambda}\\
    &\phantom{=}\left.- F (\lambda, x-L) \frac{1-e^{- A L}e^{i \lambda L}}{A- i \lambda}+2 i \lambda D \, F (\lambda, x-L) \frac{1}{D\lambda^2+B}\right] d \lambda, \\
    \frac{\partial w_2}{\partial x}  (x,t) &= -\frac{1}{2\pi}\int_{C^{+}} \frac{1}{\Delta (\lambda, -L)} \left[ 2 i \lambda D \, F (\lambda, x-L) \frac{e^{Bt}}{D\lambda^2+B}\right] d \lambda ,
\end{align*}
where 
\[
F(\lambda ,y) = \frac{\partial \Delta}{\partial y} (\lambda, y) = -\lambda^2 \sin (\lambda y) - C\lambda \cos (\lambda y).
\]

The first integral in \eqref{eq_final} can be computed efficiently numerically due to the exponential decay of the term $e^{-D\lambda^2 t}.$ On the other hand, the term $w_2$ is preferable to be computed analytically using the Cauchy or Residue theorem, depending on the appearance of singularities in the domain of integration. Different cases will be examined in the next section. 

\begin{figure}[t!]
    \centering
    \includegraphics[scale=0.8]{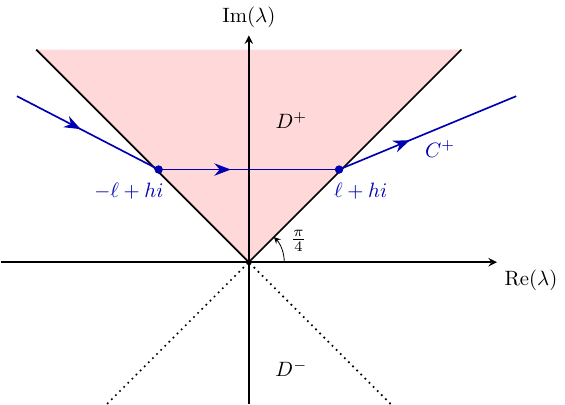}
    \caption{The boundary \(\partial D^{+}\) (black curve) of the domain \(D^{+}\) and the deformed contour  \(C^{+}\) (blue curve) given by \eqref{eq_para}.}
    \label{fig1}
\end{figure}

Initially, we choose an infinite trapezoidal contour \(C^{+}\) (positive oriented) with parametrization 
$z : \mathbb{R} \to \mathbb{C}$:
\begin{equation}\label{eq_para}
    z(s) = 
    \begin{cases} 
       z_1 (s)= (-\ell + hi) + (-s-\ell)\, e^{i\tfrac{5\pi}{6}}, & s \in (-\infty, -\ell), \\
       z_2 (s)= s + hi ,                            & s \in [-\ell,\, \ell],      \\
       z_3 (s)= (\ell + hi) + (s-\ell)\, e^{i\tfrac{\pi}{6}},    & s \in (\ell, +\infty), 
    \end{cases}
\end{equation}
as shown in Figure~\ref{fig1}.

Let $V_1(\lambda,x,t)$ denote the integrand of $w_1.$ Then, the term $w_1$ is expressed as:
\begin{align*}
\int_{C^{+}} V_1 (\lambda,x,t) d \lambda &= \int_{-\infty}^{-\ell} V_1 (z_1 (s),x,t) z_1'(s) d s +   \int_{-\ell}^\ell V_1 (z_2 (s),x,t) z_2'(s) d s  \\
&\phantom{=}+   \int_\ell^{\infty} V_1 (z_3 (s),x,t) z_3'(s) d s \\
&= \int_\ell^{\infty} \left( V_1 (z_3 (s),x,t) e^{i \tfrac{\pi}{6}} -V_1 (z_1 (-s),x,t) e^{i \tfrac{5\pi}{6}}  \right) d s + \int_{-\ell}^\ell V_1 (z_2 (s),x,t)  d s \, .
\end{align*}

\section{Numerical examples}

We verify the applicability of the analytical solution \eqref{solution_int} by comparing it with the Fourier series solution presented in \cite{clothier, hills}. The comparison reproduces the same results; therefore, we keep the same soil properties, namely  $a = 9.88\times10^{-5}\ \text{m/s},\,  b = -0.0065\ \text{m}^3/\text{m}^3,$ and $D = 3.51\times10^{-7}\ \text{m}^2/\text{s}.$

The construction of the contour $C^+$ ensures that the roots of $\Delta (\lambda,-L)$, which all lie on the real axis, are avoided. From \eqref{eq_final} we observe that the denominators vanish at $\lambda = \pm A i$ and at $\lambda =  \pm \gamma \, i,$ with $\gamma = \sqrt{\tfrac{B}{D}}.$ The negative roots are outside the domain of integration, and $\lambda = A i$ is a removable singularity since the numerator of the relevant ratio also vanishes (see the second integral in \eqref{solution_integ}).

In the first example, the finite column has length $L = 0.25\ \text{m}$, and the constant flux on the surface is given by $q=3.4 \times 10^{-6} \ \text{m/s}.$ The initial and boundary water moistures values are $\theta_0 = \theta_L = 0.03.$

\begin{figure}[t]
    \centering
    \begin{minipage}[t]{0.4\textwidth}
        \centering        \includegraphics[width=\linewidth]{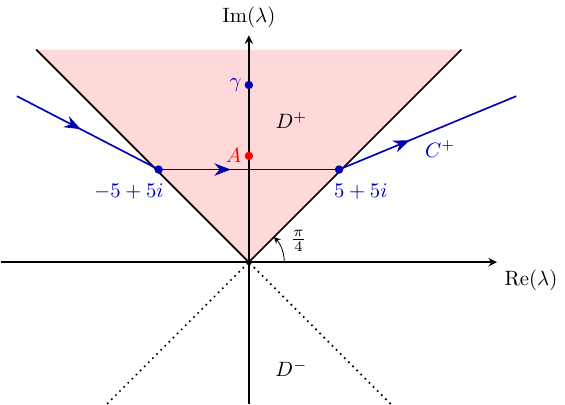}
        \subfloat{(a) Eq. \eqref{eq_para} with $\ell = 5.$}
    \end{minipage}
    \hspace{0.1\textwidth}
    \begin{minipage}[t]{0.4\textwidth}
        \centering
        \includegraphics[width=\linewidth]{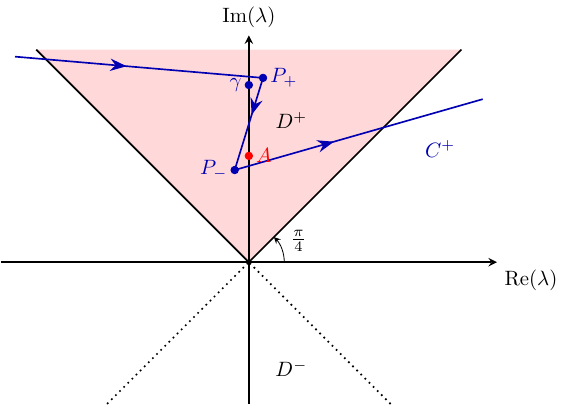}
       \subfloat{(b) Eq. \eqref{eq_para2}.}
    \end{minipage}
    \caption{The two different contour deformations $C^+$ considered in the first example. } \label{fig_ex1a}
\end{figure}

We consider the parametrization \eqref{eq_para} with $\ell=5.$ We get $A \approx 6.6 i$ and  $\gamma  \approx 52.2 i$ is a simple root in the interior of the domain with boundary $C^+$, see the left image in Figure~\ref{fig_ex1a}. Thus, we compute analytically $w_2$ using the Residue theorem, namely
\begin{align*}
w_2 (x,t) &= - \frac{1}{2\pi}\int_{C^{+}} V_2 (\lambda , x,t) d \lambda = - i \, \mbox{Res} [V_2, \gamma i] \\
&= e^{Bt} \frac{\sqrt{a q} \cosh (\gamma (L-x)) + a (\theta_L + b) \sinh (\gamma (L-x)) }{\sqrt{a q} \cosh (\gamma L) + a (\theta_L + b) \sinh (\gamma L) }.
    \end{align*}

\begin{figure}[ht]
    \centering
    \begin{minipage}[t]{0.4\textwidth}
        \centering        \includegraphics[width=\linewidth]{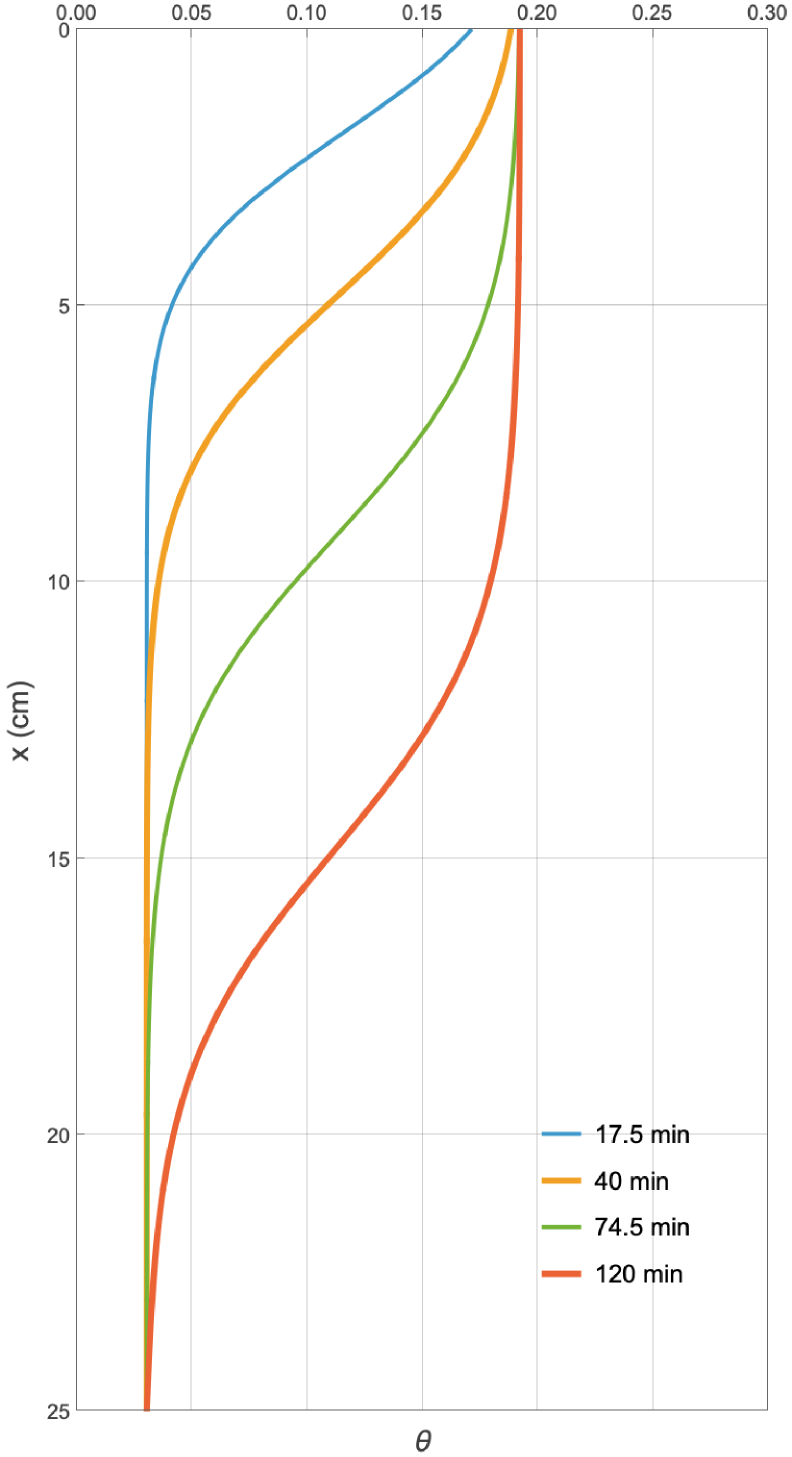}
    \end{minipage}
    \hspace{0.1\textwidth}
    \begin{minipage}[t]{0.407\textwidth}
        \centering
        \includegraphics[width=\linewidth]{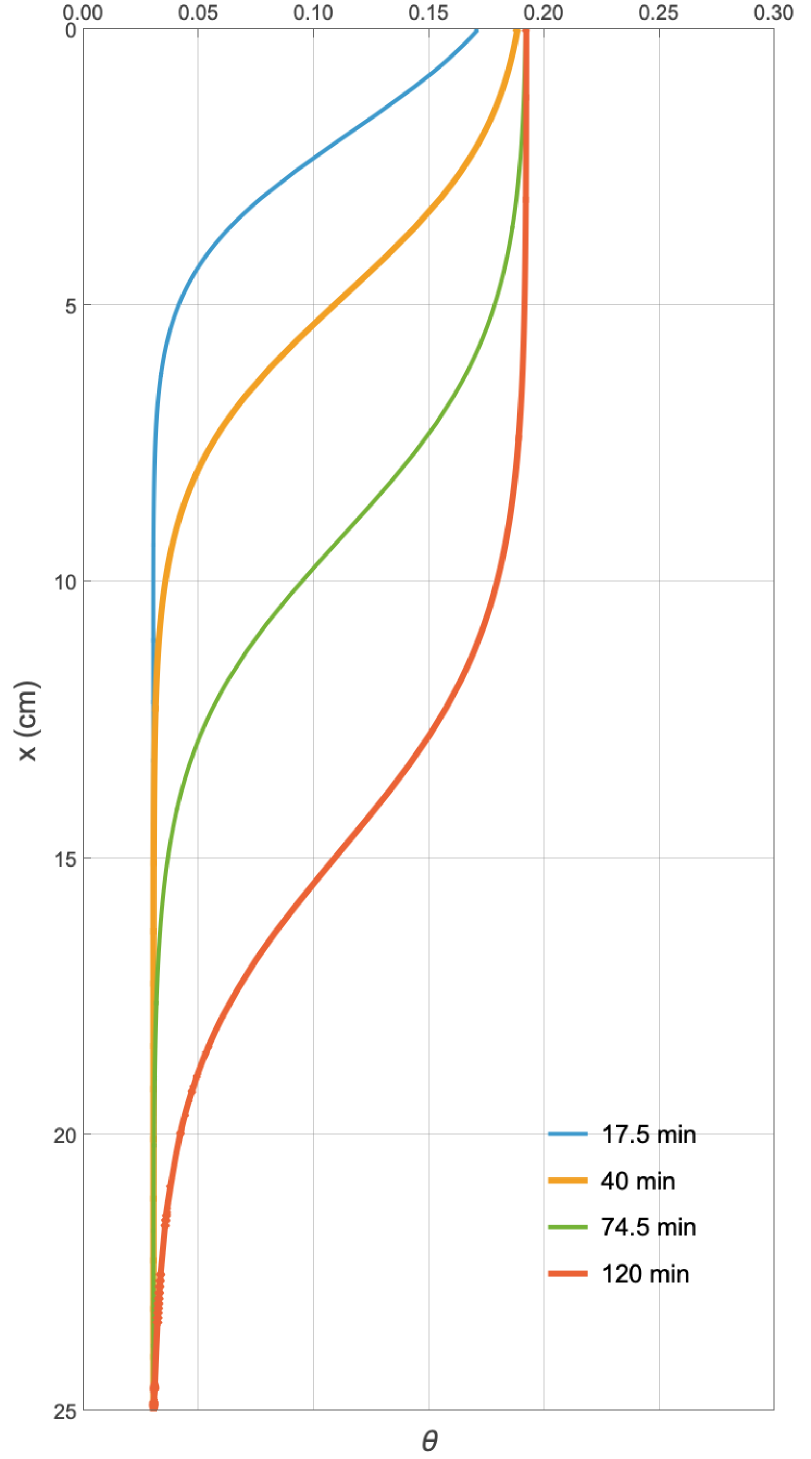}
    \end{minipage}
    \caption{The water content $\theta$ \eqref{solution} using the analytical solution \eqref{solution_int}  (left)  and the Fourier series solution, see \cite{hills}, (right) for the setup of the first example. } \label{fig_ex1_sol}
\end{figure}

In Figure~\ref{fig_ex1_sol} we present the numerical evaluation of the analytical solution, and we compare it with the Fourier solution presented in \cite{hills}. We observe that the expected physically meaningful behavior is clearly observed, and the two solutions fit perfectly.

\begin{figure}[ht]
    \centering
    \includegraphics[scale=0.45]{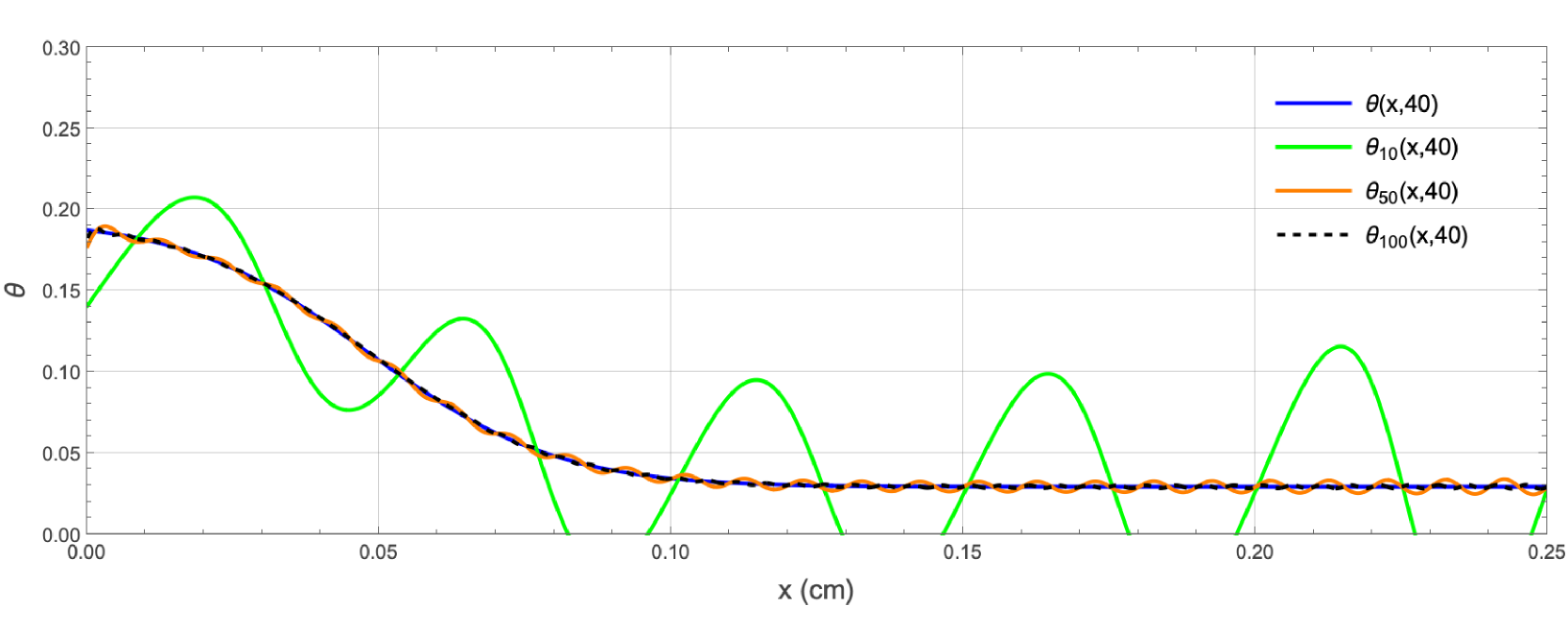}
    \caption{Convergence of the series solution to the analytical one (blue curve) for \(t=40\) min for the setup of the first example.}
    \label{fig:fouriercomp}
\end{figure}

Note that the series had to be truncated at \(N=2000\) to obtain an illustration of the solution with a sufficiently small deviation from the actual solution. A visualization of the truncation error for different numbers of Fourier nodes is shown in Figure~\ref{fig:fouriercomp} for  \(t=40\) min.

\begin{figure}[ht]
    \centering
    \includegraphics[scale=0.8]{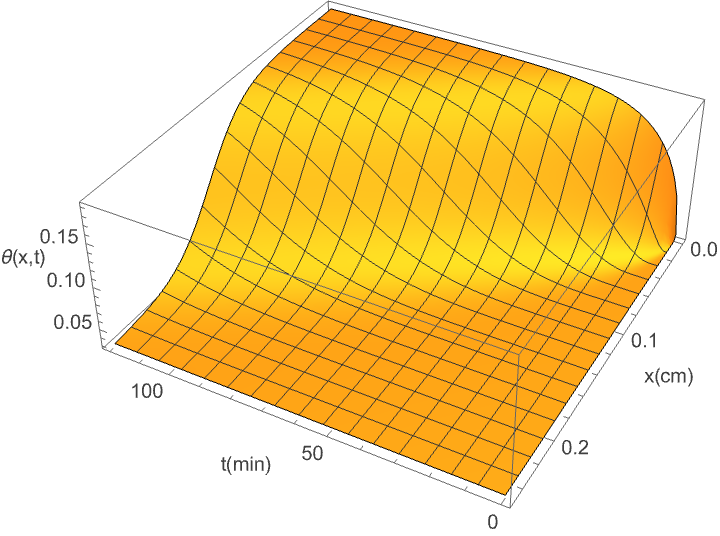}
    \caption{Example 1: The three-dimensional visualization of the analytical solution \(\theta(x,t)\) for \(x \in[0,25]\) cm and $t\in [0,120]$  min. }
    \label{fig:3dex1}
\end{figure}

For matters of broadness and illustration, in Figure ~\ref{fig:3dex1} we present the analytical solution \(\theta(x,t)\) for \(x \in[0,25]\) cm and $t\in [0,120]$ min, by adopting a different approach: We introduce an alternative parametrization of $C^+$ to avoid the computation of the  singularity at $\lambda = \gamma \, i $, namely,
\begin{equation}\label{eq_para2}
    z(s) = 
    \begin{cases} 
        P_+ + s \, e^{i\pi/8}, & s \in (-\infty,\, 0), \\
        (P_- - P_+)\,s + P_+ ,  & s \in [0,\, 1],        \\
        P_- + (s-1)\, e^{i7\pi/8}, & s \in (1,\, +\infty),
    \end{cases}
\end{equation}
where $P_- = iA - (1+i)$ and $P_+ = i \gamma+ (1+i),$ see the right image in Figure~\ref{fig_ex1a}. This way, the root \(\lambda =\gamma \, i\) lies in the exterior of the domain defined by the the boundary $C^+$, the integrand is analytic, and using Cauchy's theorem, we get
\begin{equation}\label{eq_pole}
w_2(x,t)=0.    
\end{equation}

 In the second example, the length is $L=0.08\ \text{m}$ and no-flux (outflow) is assumed, meaning $q=0.$ In this case $\theta_0 = 0.355$ and 
$\theta_L = 0.10.$ Here, $\gamma=0$ and the root is outside the region defined by \eqref{eq_para}, so \eqref{eq_pole} holds true  by Cauchy's theorem. 

\begin{figure}[ht]
    \centering
    \begin{minipage}[t]{0.4\textwidth}
        \centering        \includegraphics[width=\linewidth]{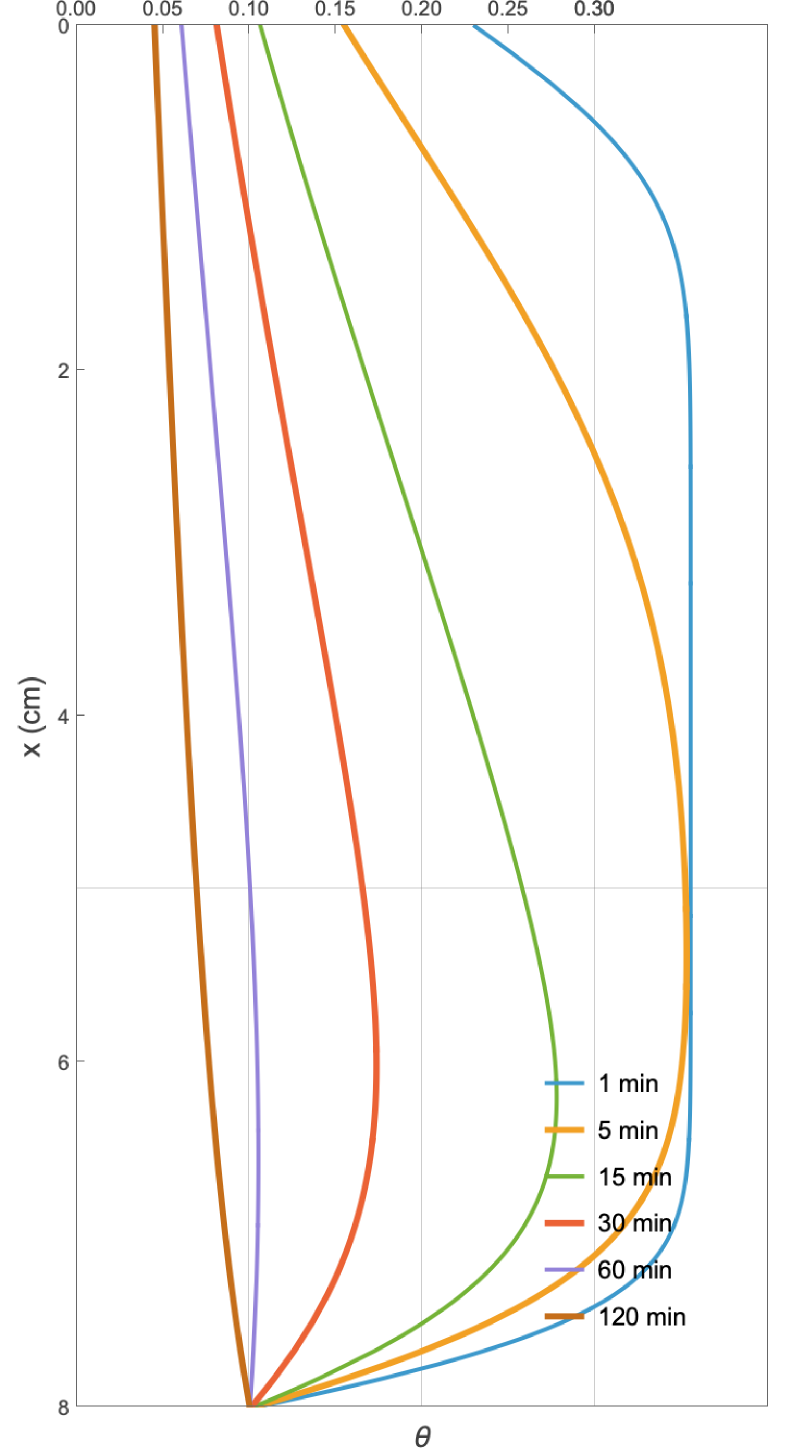}
    \end{minipage}
    \hspace{0.1\textwidth}
    \begin{minipage}[t]{0.407\textwidth}
        \centering
        \includegraphics[width=\linewidth]{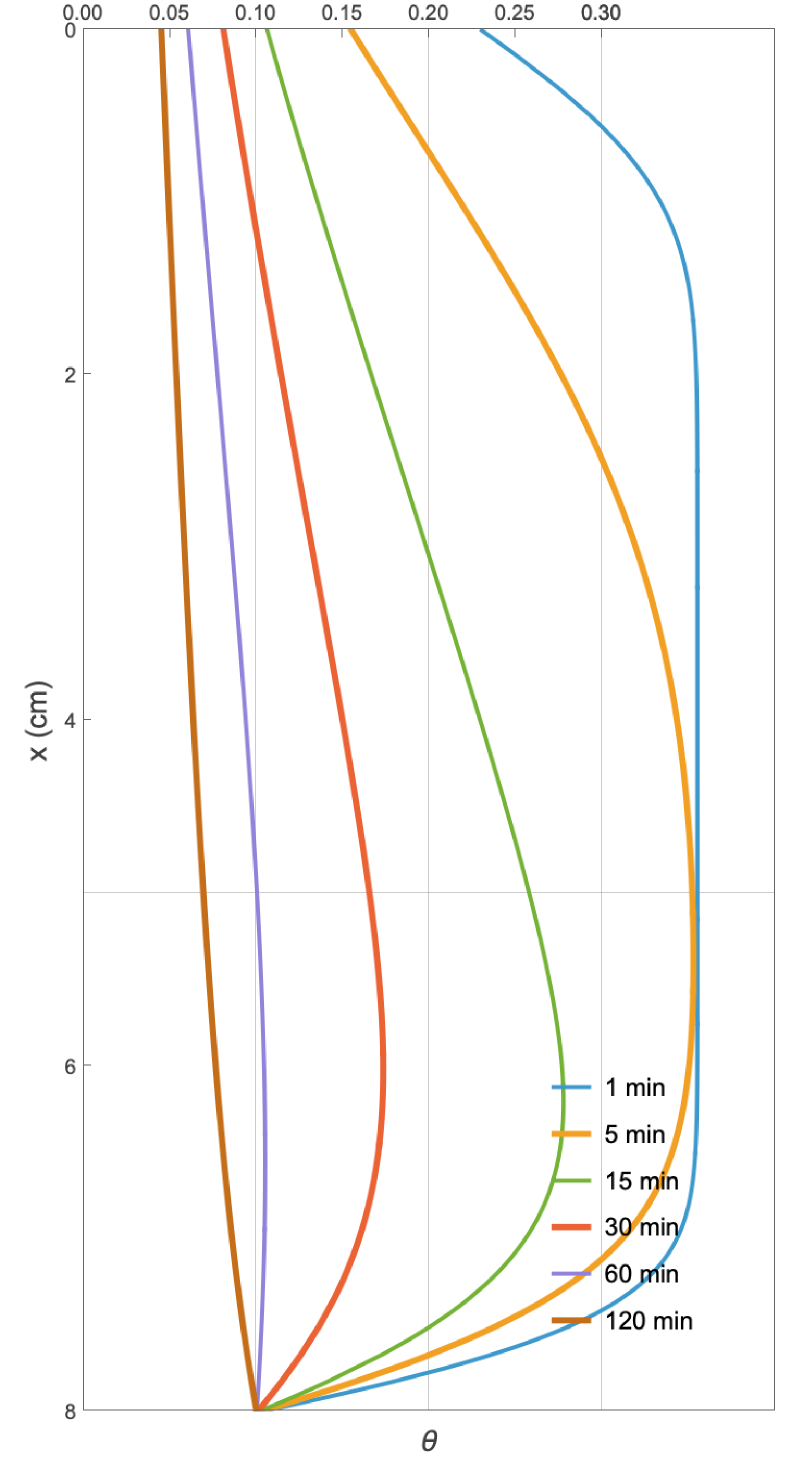}
    \end{minipage}
    \caption{The water content $\theta$ \eqref{solution} using the analytical solution \eqref{solution_int}  (left)  and the Fourier series solution, see \cite{hills}, (right) for the setup of the second example.}\label{ex2_comparison}
\end{figure}

\begin{figure}[ht]
    \centering
    \includegraphics[scale=0.8]{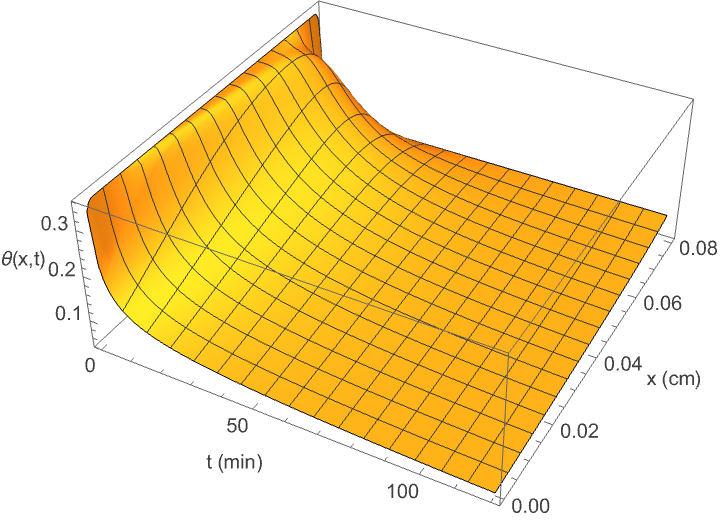}
    \caption{Example 2: The three-dimensional visualization of the analytical solution \(\theta(x,t)\) for \(x \in[0,8]\) cm and $t\in [0,120]$  min. }
    \label{fig:3dex2}
\end{figure}

The comparison between the Fokas method solution and the  Fourier series solution for \(N=2000\) is presented in Figure~\ref{ex2_comparison}.  The two solutions are again in complete agreement. A three-dimensional visualization of the analytical solution \(\theta(x,t)\) for \(x \in[0,8]\) cm and $t\in [0,120]$ min is shown in Figure~\ref{fig:3dex2}.

 \section{Conclusions}

In this work, we examined the one-dimensional vertical infiltration problem in an unsaturated porous medium modeled through the nonlinear Richards’ equation. Under certain assumptions on the soil properties, the initial boundary value problem (IBVP) was simplified to Burgers’ equation. This is still a nonlinear partial differential equation, but can be easily linearized using the Hopf–Cole transformation. This resulted in a diffusion equation on a finite interval with Robin-Dirichlet boundary conditions.

Building on previous works by the two authors, we solved the simplified IBVP using Fokas method. The key advantage is that we obtained an explicit integral representation of the solution without relying on eigenfunction expansions, with the limitations of slowly converging series. 

We tested our solution on two numerical experiments modeling real-world applications of subsurface water flow under different ponding conditions. By exploiting tools from complex analysis, we derived efficient numerical schemes that are easy to implement. The results show perfect agreement with existing solutions while offering clear practical advantages, particularly in capturing short-time behavior and near-boundary dynamics.

As future work, it would be of interest to extend this approach to more general nonlinearities, multidimensional settings, and time-dependent boundary conditions that more accurately reflect realistic hydrological scenarios. Applications to layered soil structures also appear feasible.

\bibliographystyle{abbrv}
\bibliography{refs}

\end{document}